\documentclass[letterpaper, 10 pt, conference]{ieeeconf}
\IEEEoverridecommandlockouts                              %
\usepackage[table]{xcolor}

\usepackage{graphicx}
\usepackage{epsfig}
 \usepackage{amsmath}
\usepackage{amssymb}
\usepackage{amsfonts}
\usepackage{amsthm}
 \usepackage{epstopdf}
\usepackage{algorithm}
\usepackage{algpseudocode}
\usepackage{cite}
\usepackage{subfigure}

\newcommand{\tr}{\operatorname{tr}}
\newcommand{\rank}{\operatorname{rank}}
\newcommand{\Rs}{\mathbb R}
\newcommand\rd{\cellcolor{red!20}}
\newcommand\gr{\cellcolor{green!20}}
\newcommand\bl{\cellcolor{blue!20}}
\title{\LARGE \bf
Low-rank Kalman filtering under model uncertainty
}

\author{Shenglun Yi and Mattia Zorzi \thanks{S. Yi  is with the School of Automation, Beijing Institute of Technology, Beijing 100081, China; M. Zorzi is with the Department of Information Engineering, University of Padova, Via Gradenigo 6/B, 35131 Padova, Italy. Emails: {\tt\small 3120185460@bit.edu.cn}, {\tt\small zorzimat@dei.unipd.it}}
}

\begin{document}

\maketitle
\thispagestyle{empty}
\pagestyle{empty}

\begin{abstract}  We consider a robust filtering problem where the nominal state space model is not reachable and different from the actual one. We propose a robust Kalman filter which solves a dynamic game: one player selects the least-favorable model in a given ambiguity set, while the other player designs the optimum filter for the least-favorable model. It turns out that the robust filter is governed by a low-rank risk sensitive-like Riccati equation. Finally, simulation results show the effectiveness of the proposed filter.
\end{abstract}

\section{Introduction}
Kalman filtering is widely used. However, it leads to poor performances in some applications. This is because the filter is based on a nominal model which is usually different  than the actual one. To address such a weakness, robust versions based on the standard Kalman filter have been considered, \cite{HASSIBI_SAYED_KAILATH_BOOK,speyer_2008,GHAOUI_CALAFIORE_2001,kim2020robust,Li2020}.

Risk sensitive Kalman filters \cite{RISK_WHITTLE_1980,RISK_PROP_BANAVAR_SPEIER_1998,LEVY_ZORZI_RISK_CONTRACTION,huang2018distributed}
address the model uncertainty  by replacing the standard quadratic loss function by an exponential quadratic loss function. The latter severely penalizes large errors. Such severity is tuned by the so called risk sensitivity parameter. Later on, it has been proved that risk sensitive Kalman filtering is equivalent to a minimax game, \cite{boel2002robustness,HANSEN_SARGENT_2005,YOON_2004}: one player (called nature) selects the least-favorable model in a given set (called ambiguity set), while the other player designs the optimum filter for the least-favorable model. Here, the ambiguity set is a ball, in the Kullback-Leibler (KL) divergence topology, about the nominal model. The radius of this ball depends on how much uncertainty the nominal model contains.

A modern formulation of risk sensitive filters is represented by robust Kalman filters where the model uncertainty is expressed incrementally, \cite{ROBUST_STATE_SPACE_LEVY_NIKOUKHAH_2013,zorzi2018robust,abadeh2018wasserstein,zorzi2019distributed,RKDISTR_opt,RS_MPC_IET}. More precisely, the ambiguity set is specified at each time step: in this way  the nature cannot concentrate all the uncertainty in one specific time step. These filters are built from a robust static estimation problem  showing that the Bayes estimator is optimal with respect to the ambiguity set formed by the KL divergence, \cite{robustleastsquaresestimation}. Interestingly, these results can be extended to ambiguity sets formed by the $\tau$-divergence, \cite{STATETAU_2017,OPTIMALITY_ZORZI}. These robust Kalman filters, however, are well defined
only in the case that the dynamic game involves non-degenerate probability densities. Such a condition is guaranteed by assuming that the nominal state space model is reachable and observable.

The contribution of this paper is to extend this robust Kalman filtering approach in the case that reachability is not assumed. In this case the dynamic game could involve degenerate Gaussian probability densities, i.e. their covariance matrices are low-rank. It turns out that the resulting robust Kalman filter is governed by a low-rank risk sensitive-like Riccati iteration. Although low-rank and singular Riccati iterations have been studied in the literature,
e.g. \cite{bonnabel2013geometry,Ferrante20141176,ferrante2013generalised}, our iteration appears to be new.

The outline of the paper is as follows. Section \ref{sec_2} regards the low-rank robust static estimation problem showing the optimality of the Bayes estimator. The latter is then used to derive the low-rank robust Kalman filter in Section \ref{section_3}. Some simulations are presented in Section \ref{sec_4} which show that the proposed filter outperforms the Kalman filter in the case of model uncertainty. Finally, the conclusions and the future work are discussed in Section \ref{sec_5}.

We warn the reader that the present paper only reports some preliminary result regarding the robust estimation under model uncertainty in the case that reachability is not assumed. In particular, all the proofs and most of the technical assumptions needed therein are omitted and will be published afterwards.

{\em Notation:}  The image of matrix $K$ is denoted by $\mathrm{Im}(K)$. Given a symmetric matrix $K$:  $K>0$ ($K\geq 0$) means that $K$ is positive (semi) definite; $\sigma_{max}(K)$ is the maximum eigenvalue of $K$. The symbol $\operatorname{diag}\left(d_1, \cdots, d_{n} \right)$ denotes the diagonal matrix whose entries in the main diagonal are $d_1\dots d_n$. $x\sim \mathcal N(m,K)$ means that $x$ is a Gaussian random variable with mean $m$ and covariance matrix $K$.
\section{ Low-rank robust static estimation}\label{sec_2}

We consider a robust static estimation problem where we want to estimate a random vector $x$, taking values in $ \mathbb{R}^{n}$, given the observation vector $y$, taking values in $\mathbb{R}^{p}$, and whose joint probability distribution is degenerate. More precisely, let $z:=\left[\begin{array}{ll}{x^{T}} & {y^{T}}\end{array}\right]^{T}$ and $f(z) \sim \mathcal{N}\left(m_{z}, K_{z}\right)$ denote the nominal probability density function of $z$ where $m_z\in\mathbb R^{n+p}$ and $K_z\in\mathbb R^{n+p\times n+p}$ are such that
$$
m_{z}=\left[\begin{array}{c}{m_{x}} \\ {m_{y}}\end{array}\right], \quad K_{z}=\left[\begin{array}{cc}{K_{x}} & {K_{x y}} \\ {K_{y x}} & {K_{y}}\end{array}\right].
$$
We assume that $K_z$ is such that $\mathrm{rank}(K_z)=r+p$ with $r<n$ and $K_y>0$. Moreover, let $\tilde f(z) \sim \mathcal{N}(\tilde m_{z}, \tilde K_{z})$ be the actual probability density function where $\mathrm{rank}(\tilde K_z)=r+p$. Accordingly,
\begin{equation} \label{pdf_f}
\begin{aligned} f(z) =&\left[(2 \pi)^{r+p} \operatorname{det}^{+}\left(K_{z}\right)\right]^{-1 / 2}  \\
& \times  \exp \left[-\frac{1}{2}\left(z-m_{z}\right)^{T} K_{z}^{+}\left(z-m_{z}\right)\right] \end{aligned}
\end{equation}
and
\begin{equation}\label{pdf_tildef}
\begin{aligned} \tilde f(z) =&\left[(2 \pi)^{r+p} \operatorname{det}^{+}\left(\tilde K_{z}\right)\right]^{-1 / 2}   \\
& \times \exp \left[-\frac{1}{2}\left(z-\tilde m_{z}\right)^{T} \tilde K_{z}^{+}\left(z-\tilde m_{z}\right)\right] \end{aligned}
\end{equation}
where $K^{+}_z$ and $\tilde K^{+}_z$ are the pseudo-inverse of $K_z$ and $\tilde K_z$, respectively, and $\operatorname{det}^{+}$ is the pseudo-determinant. Notice that the supports of $f_z$ and $\tilde f_z$ are, respectively, the $r+p$-dimensional affine subspaces
\begin{align*}
\mathcal{A}&=\left\{m_{z}+v, \quad v \in \mathrm {Im} \left(K_{z}\right)\right\}\\
\mathcal{\tilde A}&= \{\tilde m_{z}+v, \quad v \in \mathrm {Im}(\tilde{ K}_{z}) \}.
\end{align*}
To measure the mismatch between $f(z)$ and $\tilde f(z)$ we introduce the KL-divergence between these degenerate probability density functions. In order to do that, we have to impose $\mathcal{A}=\mathcal{\tilde A}$, indeed such divergence is not able to measure ``deterministic'' discrepancies between the nominal and the actual model. The latter condition is equivalent to impose that  $$ \mathrm {Im}\left(K_{z}\right)=\mathrm {Im} (\tilde{ K}_{z}), \; \; \Delta m_z \in \mathrm {Im}\left(K_{z}\right) $$ where $\Delta m_z=\tilde m_z-m_z$. Under the aforementioned assumption, the KL-divergence is defined as
\begin{equation} \label{def_DL_ddeg}
D (\tilde{f}, f )=\int_\mathcal{A} \ln \left(\frac{\tilde{f}(z)}{f(z)}\right) \tilde{f}(z) d z.
\end{equation}
Substituting (\ref{pdf_f}) and (\ref{pdf_tildef}) in (\ref{def_DL_ddeg}), it is not difficult to see that
\begin{equation} \label{D}
\begin{aligned} D(\tilde f, {f})=& \frac{1}{2}\left[\Delta m^{T}_z K_{z}^{+} \Delta m_z+\ln \operatorname{det}^{+} (K_{z})\right.\\ &\left.-\ln \operatorname{det}^{+} (\tilde{K}_{z})+\tr\left(K_{z}^{+} \tilde{K}_{z}\right)-(r+p)\right] .\end{aligned}
\end{equation}

\newtheorem{lemma}{Lemma}
\begin{lemma}
\label{lemma1}
Let $f(z) \sim \mathcal{N}\left(m_{z}, K_{z}\right)$ and $\tilde f(z) \sim \mathcal{N}(\tilde m_{z}, \tilde K_{z})$ be degenerate Gaussian probability density functions with the same $r+p$-dimensional support $\mathcal{A}$. Let \begin{align*}\mathcal U&=\{\tilde m_z \in\mathbb R^{n+p} \hbox{ s.t. } \tilde m_z-m_z\in \mathrm{Im}(K_z)\}\\
\mathcal{V}&=\{\tilde K_{z} \in \mathbb{R}^{n+p \times n+p} ~\text {s.t.}~ K_z=K_z^{T}, \operatorname{Im}(K_{z} )=\operatorname{Im} (\tilde{K}_{z} )\}. \end{align*} Then, $D(\tilde f,f)$ is strictly convex with respect to $\tilde m_z\in\mathcal U$ and $\tilde K_z\in\mathcal V$. Moreover, $D(\tilde f, f) \geq 0$ and equality hold if and only if $f=\tilde f$.
\end{lemma}

Assume that the nominal density $f$ is known while the actual one $\tilde f $,  having the same support of $f$, is not. Next, we design a robust estimator $\hat x=g^0(y)$ of $x$ according to the worst probability density $\tilde f(z)$ belonging to the ambiguity set which is a ball
 $$
\mathcal{B}=\{\tilde{f} \sim \mathcal N(\tilde m_z,\tilde K_z)\hbox{ s.t. } D(\tilde{f}, f) \leq c\}
$$ where $c$ is the mismodeling budget hereafter called tolerance.
More precisely, we aim to solve the following minimax problem
\begin{equation}
(\tilde f^0, g^0)=\operatorname{arg}\min _{g \in \mathcal{G}} \max _{\tilde{f} \in \mathcal{B}} J(\tilde{f}, g)
\end{equation}
where$$
\begin{aligned} J(\tilde{f}, g) &=\frac{1}{2} E_{\tilde{f}}\left[\|H(x-g(y))\|^{2}\right] \\ &=\frac{1}{2} \int_\mathcal{A}\|H(x-g(y))\|^{2} \tilde{f}(z) d z; \end{aligned}
$$  $H \in \mathbb{R}^{q \times n}$ with $q \leq r$ and full row rank; $\mathcal{G}$ is the set of estimators for which $E_{\tilde{f}}\left[\|H(x-g(y))\|^{2}\right] $ is bounded with respect to all the  Gaussian densities in $\mathcal B$.

\newtheorem{theorem}{Theorem}
\begin{theorem} \label{theo1} Let $f$ be a Gaussian (possibly degenerate) density defined as in (\ref{pdf_f}) with $K_y>0$. Assume that $\mathrm {Im}(H^T) \subseteq \mathrm {Im} (P)$ with
\begin{align} \label{def_P_stat}P&:=K_{x}-K_{x y} K_{y}^{-1} K_{y x}.\end{align} Then, the least favorable Gaussian density $\tilde f^0$ is with mean vector and covariance matrix
$$
\tilde{m}_{z}^0=m_{z}=\left[\begin{array}{c}{m_{x}} \\ {m_{y}}\end{array}\right], \quad \tilde{K}_{z}^0=\left[\begin{array}{cc}{\tilde{K}_{x}} & {K_{x y}} \\ {K_{y x}} & {K_{y}}\end{array}\right]
$$
so that, only the covariance of $x$ is perturbed.
Then, the optimal robust estimator is the Bayes estimator
\begin{equation}
g^{0}(y)=G_{0}\left(y-m_{y}\right)+m_{x}
\end{equation}
with $G_0=K_{x y} K_{y}^{-1}$.
 The nominal posterior covariance matrix of $x$ given $y$ is given by (\ref{def_P_stat}), while the least favorable one is:
\begin{align*}
 \tilde P &:=\tilde{K}_{x}-K_{x y} K_{y}^{-1} K_{y x}.
\end{align*}  Then, we have $$\tilde P=(P^+-\lambda ^{-1}H^T H)^+.$$
Moreover, there exists a unique Lagrangian multiplier $\lambda >\sigma_{max}(Q)>0$, such that $c=D(\tilde f^0, f)$ where
 $$Q:=HH^T(HP^+H^T)^{-1}HH^T.
$$
\end{theorem}

{\em Remark:} In the case that $P>0$ then $f$ is a non-degenerate density. In such a case, Theorem \ref{theo1} still holds and: the pseudo-inverse is replaced by the inverse; moreover, $H^T H$ can be chosen as the identity matrix. In the latter case, we recover the robust static estimation problem proposed in \cite{robustleastsquaresestimation}.

\section{ Low-rank robust Kalman filter} \label{section_3}
We consider a nominal Gauss-Markov state space model of the form:
\begin{equation}\label{ss_space}
\left\{\begin{array}{c}{x_{t+1}=A_t x_t+B_t v_t} \\ {y_t=C_t x_t+D_t v_t}\end{array}\right.
\end{equation}
where $A_t\in\Rs^{n\times n}$, $B_t\in\Rs^{n\times n+p}$, $C_t\in\Rs^{p\times n}$, $D_t\in\Rs^{p\times n+p}$; $x_t$ and $y_t $ are the state vector and the observation vector, respectively. $v_t$ is normalized white Gaussian noise. We assume that $D_tD_t^T>0$ that is all the components of the observation process are affected by a full rank $p$-dimensional noise. Let $z_t:=\left[\begin{array}{ll}{x_{t+1}^{T}} & {y_t^{T}}\end{array}\right]^{T},$ so the nominal conditional transition probability density function of the nominal state space model is $\phi_t(z_t|x_t) \sim \mathcal{N}\left(m_{z_t|x_t}, K_{z_t|x_t}\right)$ with$$
m_{z_t|x_t}=\left[\begin{array}{c}{ A_t} \\ { C_t}\end{array}\right] x_{t}, ~~~K_{z_t|x_t}=\left[\begin{array}{cc}{ B_t  B_t^T} & {B_t D_t^T} \\ { D_tB_t ^T} & { D_t  D^T_t}\end{array}\right].$$
Notice that $K_{z_t|x_t}$ is not necessarily positive, i.e. $\phi_t(z_t|x_t)$ could be degenerate. Indeed no assumption on $B_t$ has been made. Then, we assume at time $t$ the a priori conditional density of $x_t$ given $ Y_{t-1}:=\{y_0 \dots y_t\}$ is
\begin{align}\label{hp_ft}
\check{f}_{t}\left(x_{t} | Y_{t-1}\right) \sim \mathcal{N} (\hat{x}_{t}, \tilde P_{t})
\end{align}
with $\rank(\tilde P_t)=r_{t}$. Therefore, we obtain the pseudo-nominal density
\begin{equation*}
{ {f}_{t}\left(z_{t} | Y_{t-1}\right)}=\int_{\check{\mathcal A}_t} \phi_{t}\left(z_{t} | x_{t}\right){ \check{f}_{t}\left(x_{t} | Y_{t-1}\right)} d x_{t}
\end{equation*} where ${\check{\mathcal A}}_t$ denotes the support of $\check f_t$. Then, it is not difficult to see that  ${ {f}_{t}\left(z_{t} | Y_{t-1}\right)} \sim \mathcal{N}\left(m_{z_t}, K_{z_{t}}\right)$ with
\begin{equation*}
m_{z_t}=\left[\begin{array}{l}
A_{t} \\
C_{t}
\end{array}\right] \hat{x}_{t}, \quad {K}_{z_{t}}=\left[\begin{array}{cc}
{K}_{x_{t+1}} & K_{x_{t+1} y_{t}} \\
K_{y_{t} x_{t+1}} & K_{y_{t}}
\end{array}\right]
\end{equation*}
where the conditional covariance matrix $K_{z_t}$ takes the parametric form
\begin{equation}\label{K_z}
K_{z_{t}}=\left[\begin{array}{c}
A_{t} \\
C_{t}
\end{array}\right] \tilde P_{t}\left[\begin{array}{cc}
A_{t}^{T} & C_{t}^{T}
\end{array}\right]+\left[\begin{array}{c}
B_{t} \\
D_{t}
\end{array}\right]\left[\begin{array}{cc}
B_{t}^{T} & D_{t}^{T}
\end{array}\right].
\end{equation}
Notice that the support of ${f}_{t}\left(z_{t} | Y_{t-1}\right)$ is the affine subspace $$\mathcal{A}_t:=\left\{m_{z_t}+{ w_t}, \quad { w_t} \in \operatorname{Im}\left(K_{z_t}\right)\right\}.$$

Let $\tilde \phi_t(z_t|x_t)$ be the actual least favorable density of $z_t$ given $x_t$. Then, the marginal density is
\begin{align*}
\tilde{f}_{t}\left(z_{t} | Y_{t-1}\right)=\int_{{\check{\mathcal A}}_t}\tilde{\phi}_{t}\left(z_{t} | x_{t}\right) {\check{f}_{t}\left(x_{t} | Y_{t-1}\right)} d x_{t}.
\end{align*}
In what follows, we assume $\tilde \phi_t$ is Gaussian, accordingly the pseudo-actual density $\tilde{f}_{t}(z_{t} | Y_{t-1}) \sim \mathcal{N}(\tilde m_{z_t}, \tilde{K}_{z_t})$ is Gaussian. In order to measure the mismatch between $ f_t(z_t|Y_{t-1})$ and $\tilde f_t(z_t|Y_{t-1})$ using the KL-divergence we impose that $\tilde f_t(z_t|Y_{t-1})$ has the same support of  $f_t(z_t|Y_{t-1})$. Accordingly, $\mathrm {Im}(\tilde K_{z_t})=\mathrm {Im}(K_{z_t})$ and $\Delta m_{z_t}:=\tilde m_{z_t}-m_{z_t}\in\mathrm {Im}(K_{z_t})$. Under the above assumptions, we assume that the actual density belongs the following ambiguity set
\begin{align}\label{def_Bt}
\mathcal{\tilde B}_t=\{ \tilde f_t\sim \mathcal N(\tilde m_{z_t},\tilde K_{z_t}) ~~  s.t.\;\;   D(\tilde f_t,  f_t) \leq c_t\}.
\end{align}
It is worth noting that the model uncertainty in (\ref{def_Bt}) is expressed incrementally. In plain words, the tolerance $c_t$ represents the mismodeling budget allowed at time $t$. Then, we consider as robust one step-ahead predictor $\hat x_{t+1}$ of $x_{t+1}$ given $Y_t$, the solution to the following minimax game
\begin{equation} \label{minimax_dyn}
\hat x_{t+1} =\underset{g_{t} \in \mathcal{G}_{t}}{\mathrm{argmin}} \max _{\tilde{f}_{t} \in \mathcal{B}_{t}} \bar J_{t}(\tilde f_{t}, g_{t} )
\end{equation}
where $$
\bar{J}_{t}(\tilde{f}_{t}, g_{t})=\int_{\mathcal A_t}\left\|H_t\left(x_{t+1}-g_{t}\left(y_{t}\right)\right)\right\|^{2} \tilde{f}_{t}\left(z_{t} | Y_{t-1}\right) \mathrm{d} z_{t},
$$ $H_t\in \mathbb R^{q\times n}$ with $q \leq r $ is a square root of the projection matrix having the same image of  $P_{t+1}$, i.e. $H^T_tH_t$ is the projection matrix such that $\mathrm {Im} (H_t^TH_t)= \mathrm {Im} (P_{t+1})$, where
\begin{align}\label{Pt1}
P_{t+1} &:=K_{x_{t+1}}-K_{x_{t+1}, y_{t}} K_{y_{t}}^{-1} K_{y_{t}, x_{t+1}}.
\end{align}
In the case that $P_{t+1}>0$, then $H_t$ is an orthogonal matrix and thus $\|H_t(x_{t+1}-g_t(y_t))\|^2=\|(x_{t+1}-g_t(y_t))\|^2$. Accordingly, the minimax problem in (\ref{minimax_dyn}) boils down to the one in \cite{ROBUST_STATE_SPACE_LEVY_NIKOUKHAH_2013}.

{\em Remark:}  Condition (\ref{hp_ft}) means that, using the terminology coined by Hansen and Sargent \cite{ROBUSTNESS_HANSENSARGENT_2008}, the maximizer in (\ref{minimax_dyn}) is operating under commitment, i.e. the maximizer is required to commit all the least favorable model components at early stages  with the estimating player.

To solve Problem (\ref{minimax_dyn}) we use Theorem \ref{theo1}. Indeed, replacing $f$, $\tilde f$, $g$, $H$ by $f_t$, $\tilde f_t$, $g_t$, $H_t$, respectively, it is not difficult to see that all the assumptions are satisfied. In particular the condition on $H_t$ is satisfied and  we have
$$ K_{y_t}= C_t\tilde P_t C_t^T+D_tD_t^T\geq D_tD_t^T>0.$$ Then, the least favorable density is $\tilde{f}_{t}^{0}(z_{t} | Y_{t-1}) \sim \mathcal{N}(m_{z_t}, \tilde{K}_{z_t})$ where $$
\tilde{K}_{z_{t}}=\left[\begin{array}{cc}
\tilde{K}_{x_{t+1}} & K_{x_{t+1}, y_{t}} \\
K_{y_{t}, x_{t+1}} & K_{y_{t}}
\end{array}\right].
$$
Then, the nominal posterior covariance of $x_{t+1}$ given $Y_{t}$ has been defined in (\ref{Pt1}) and the the least favorable one is
\begin{align*}
\tilde P_{t+1} &=\tilde{K}_{x_{t+1}}-K_{x_{t+1}, y_{t}} K_{y_{t}}^{-1} K_{y_{t}, x_{t+1}}.
\end{align*}
Moreover, \begin{equation*}
\tilde P_{t+1}=\left(P_{t+1}^{+}-\lambda_{t+1}^{-1} H^{T}_t H_t\right)^{+}
\end{equation*}
where $\lambda_t>\sigma_{max}(H_tH_t^T(H_tP_{t+1}^{+}H_t^T)^{-1}H_tH_t^T)$ is the unique solution of the following equation
\begin{equation*}
 \begin{aligned}
  \gamma_t(P_{t+1},\lambda_t):=&\frac{1}{2}\left\{\operatorname{lndet}^+(P_{t+1})-\operatorname{lndet}^+(\tilde P_{t+1})\right.\\
&\left.+\tr\left[P_{t+1}^+\tilde P_{t+1}-I_{r_{t+1}}\right]\right\}=c_t.
\end{aligned}
\end{equation*}
Moreover, the robust estimator is
\begin{equation*}
\hat{x}_{t+1}=A_{t} \hat{x}_{t}+G_{t}\left(y_{t}-C_{t} \hat{x}_{t}\right).
\end{equation*}
where $G_t=K_{x_{t+1},y_t}K^{-1}_{y_t}$ and by equation (\ref{K_z}), we get the parametric form of $G_t$ and $P_{t+1}$ \begin{equation*}
G_{t}= (A_{t} \tilde P_{t} C_{t}^{T}+B_{t} D_{t}^{T} )(C_{t} \tilde P_{t} C_{t}^{T}+D_{t} D_{t}^{T})^{-1}
\end{equation*}
\begin{equation}\label{ric_dist}
P_{t+1}=A_{t} \tilde P_{t} A_{t}^{T}-G_{t} (C_{t} \tilde P_{t} C_{t}^{T}+D_{t} D_{t}^{T}) G_{t}^{T}+B_{t} B_{t}^{T}.
\end{equation}

\begin{algorithm}[h]
  \caption{Low-rank robust Kalman filter at time $t$}
  \begin{algorithmic}[1]
    \Require
      $y_t$, $\hat x_t$,      $\tilde P_t$, $c_t$
      \State $ G_{t}=(A_{t} \tilde P_{t} C_{t}+B_{t} D_{t}^{T})^{T}(C_{t} \tilde P_{t} C_{t}^{T}+D_{t} D_{t}^{T})^{-1}$
      \State $\hat{x}_{t+1}=A_{t} \hat{x}_{t}+G_{t}\left(y_{t}-C_{t} \hat{x}_{t}\right)$
       \State \label{RIcc_step}$P_{t+1}=A_{t} \tilde P_{t} A_{t}^{T}-G_{t}\left(C_{t} \tilde P_{t} C_{t}^{T}+D_{t} D_{t}^{T}\right) G_{t}^{T}+B_{t} B_{t}^{T}$
      \State Select $H^T_tH_t$ as projection matrix with image $ \mathrm {Im} (P_{t+1})$
      \State Find $\lambda_t$ s.t. $ \gamma_t( P_{t+1},\lambda_t)=c_t$
      \State \label{RIcc_step2} $\tilde P_{t+1}=\left(P_{t+1}^{+}-\lambda_{t}^{-1} H_t^{T} H_t\right)^{+}$
  \end{algorithmic}\label{code:recentEnd}
\end{algorithm}
Algorithm \ref{code:recentEnd} summarizes the robust Kalman filter that we obtain. It is worth noting that steps \ref{RIcc_step} and \ref{RIcc_step2} represent a risk sensitive-like Riccati iteration which is well defined also in the case that $\tilde P_t$ is low-rank. Notice that $\theta_t:=\lambda_t^{-1}$ represents the time-varying risk sensitivity parameter. In the situation that $c_t=0$, i.e. the actual model corresponds to the nominal model, then we have $\theta_t=0$ and thus $P_t=\tilde P_t$; in particular, (\ref{ric_dist}) become the  usual Riccati equation and thus we obtain the standard Kalman filter.

\section{Simulation results}\label{sec_4}
We consider the linear time-invariant model
\begin{equation}\label{m_nomi}
\left\{\begin{array}{c}{x_{t+1}=A x_t+B v_t} \\ {y_t=C x_t+D v_t}\end{array}\right.
\end{equation}
where
\begin{align}
A&=\left[\begin{array}{ccc}
1.1 & 0.1 & 0.1 \\
0 & 0.6364 & -0.6364\\
0 & 0.6364 & 0.6364
\end{array}\right],\; \;  B=\left[\begin{array}{cc}
1 & 0\\
0 & 0\\
0 & 0
\end{array}\right], \nonumber\\
C&=\left[\begin{array}{lll}
1 & 0 & 0
\end{array}\right], \; \;  D=\left[\begin{array}{ll}
 0 & 1
\end{array}\right], \nonumber
\end{align} and $x_{0} \sim \mathcal{N}\left(0, \tilde P_{0}\right)$ with
$$ \tilde P_0=\left[\begin{array}{lll}
1 & 0 & 0 \\
0 & 0 & 0 \\
0 & 0 & 1
\end{array}\right].$$ Notice that the pair $(A,B)$ is stabilizable, but not reachable. Finally, the pair $(A,C)$ is observable.

First, we compare $P_t$ and $\tilde P_t$ by using the robust Kalman filter (RKF) of section \ref{section_3} and the Kalman filter (KF). More precisely, we consider two values of the tolerance for RKF: $c_1=5 \cdot 10^{-2}$,  $c_2=8 \cdot 10^{-2}$ and we denote the corresponding robust filters as RKF1 and RKF2, respectively. \begin{figure}[h]
\centering
\includegraphics[width=0.5\textwidth]{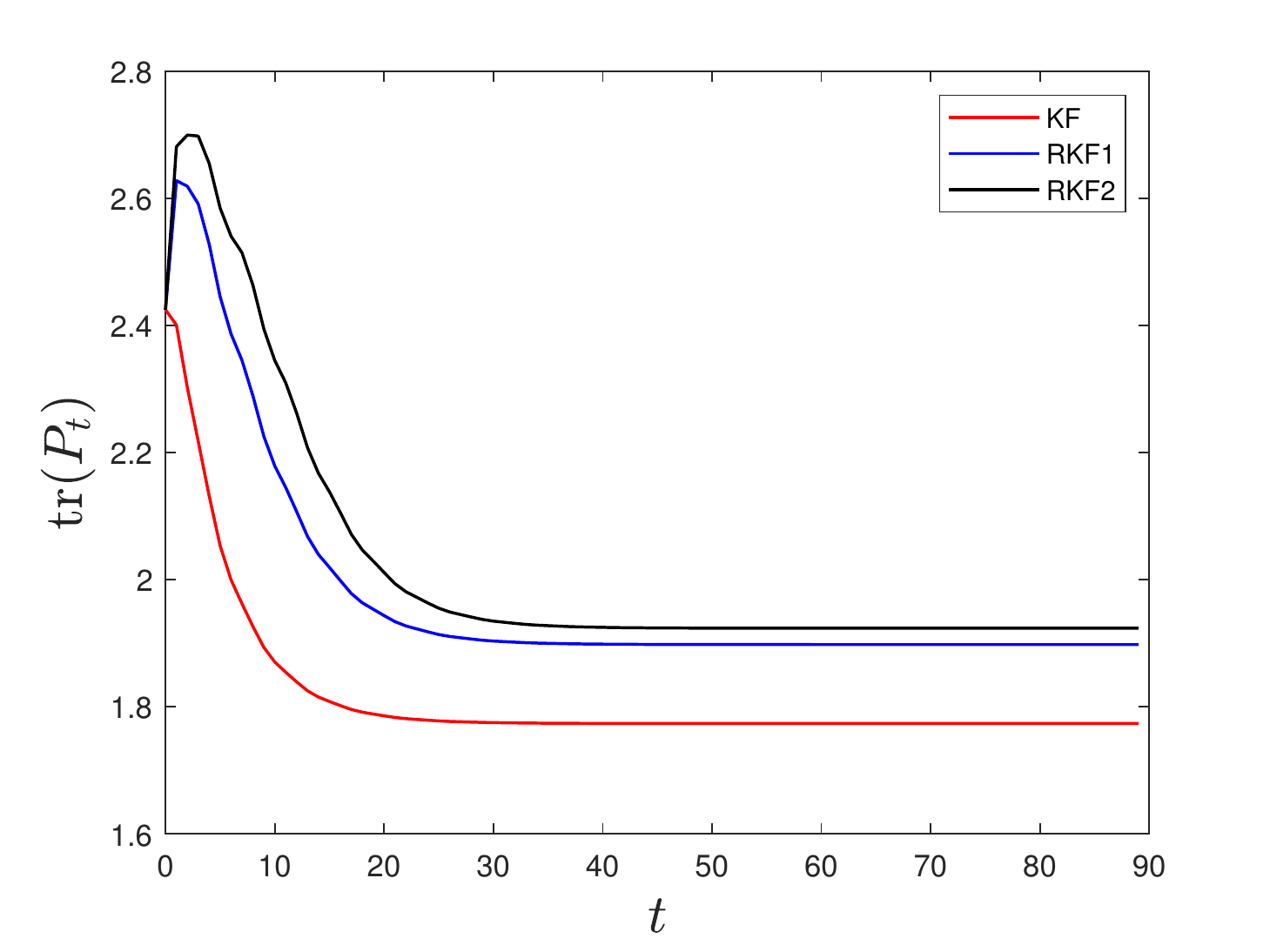}
\caption{Trace of $P_t$ for KF, RKF1 with $c=5 \cdot 10^{-2}$, and RKF2 with $c=8 \cdot 10^{-2}$. Recall that $\tilde P_t=P_t$ for KF.} \label{fig4}
\end{figure} Figure \ref{fig4} shows the trace of $P_t$ over the time horizon $[0,90]$. Clearly $\tr( P_t)$ of KF converges to a constant value because the nominal model is stabilizable and observable, indeed it is well known that the corresponding Riccati equation converges to a unique solution. In regard to RKF1 and RKF2, $\tr( P_t)$ converges for both. Moreover, the larger $c$ is, the more $\tr(P_t)$ is different from the one of KF. It is worth noticing that $\rank(P_t)=2$ for KF, RKF1 and RKF2.  \begin{figure}[h]
\centering
\includegraphics[width=0.5\textwidth]{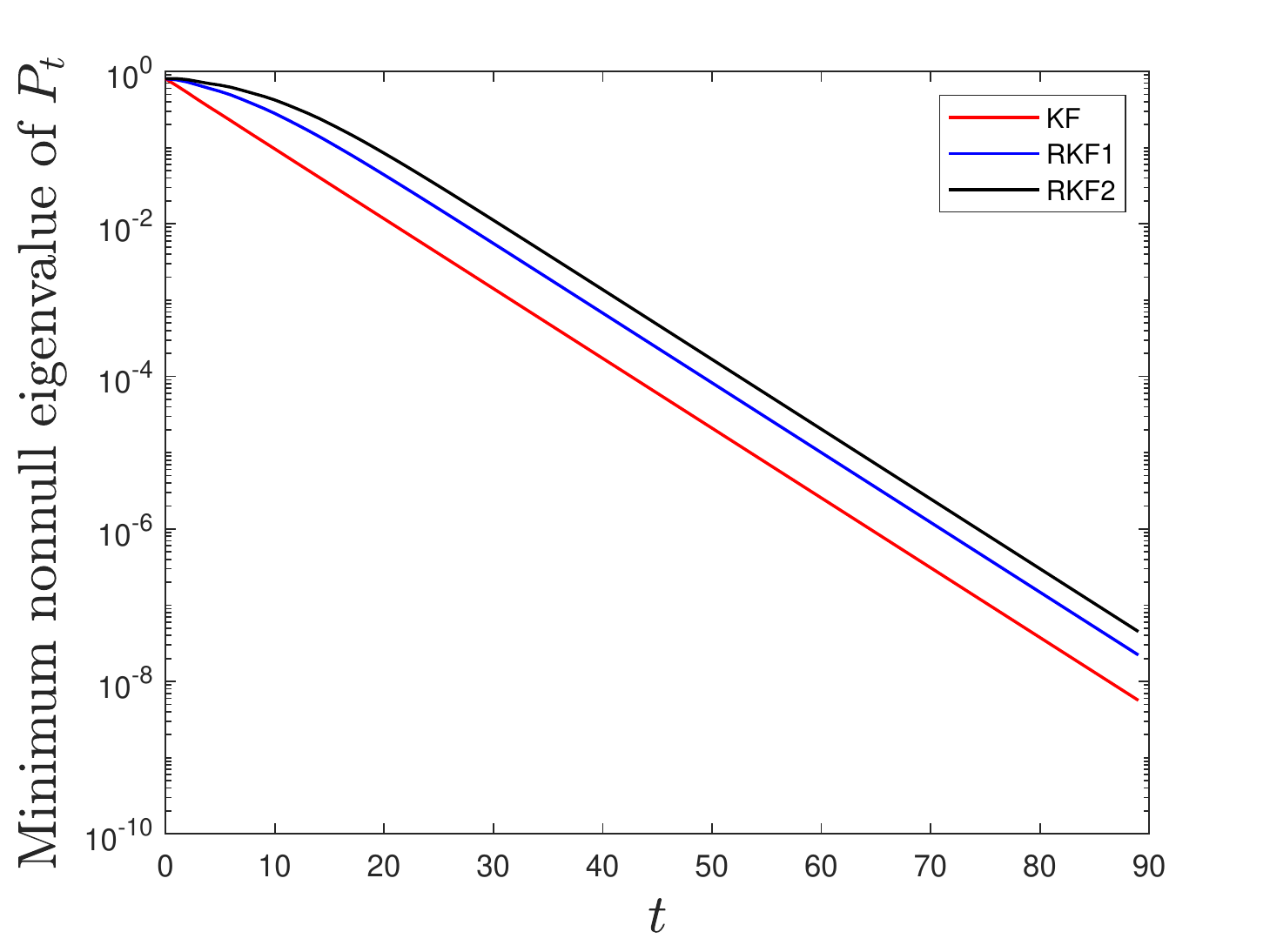}
\caption{Minimum nonnull eigenvalue of $P_t$ as a function of $t$ for KF, RKF1 with $c=5 \cdot10^{-2}$, and RKF2 with $c=8 \cdot 10^{-2}$.} \label{fig5}
\end{figure}
Figure \ref{fig5} shows the minimum nonnull eigenvalue of $ P_t$: we notice that $\rank(P_t)\rightarrow 1$ as $t\rightarrow \infty$ for all the filters. This is clearly expected from KF because the corresponding algebraic Riccati equation admits a unique solution with rank equal to one. It is also worth noting that the highest convergence rate to the rank one solution is given by KF, while the larger $c$ is, the slower the convergence rate is. Figures \ref{fig1} and \ref{fig2} show the trace of $\tilde P_t$ and its minimum nonnull eigenvalue as a function of $t$. Similar observations can be made also in this case.  \begin{figure}[h]
\centering
\includegraphics[width=0.5\textwidth]{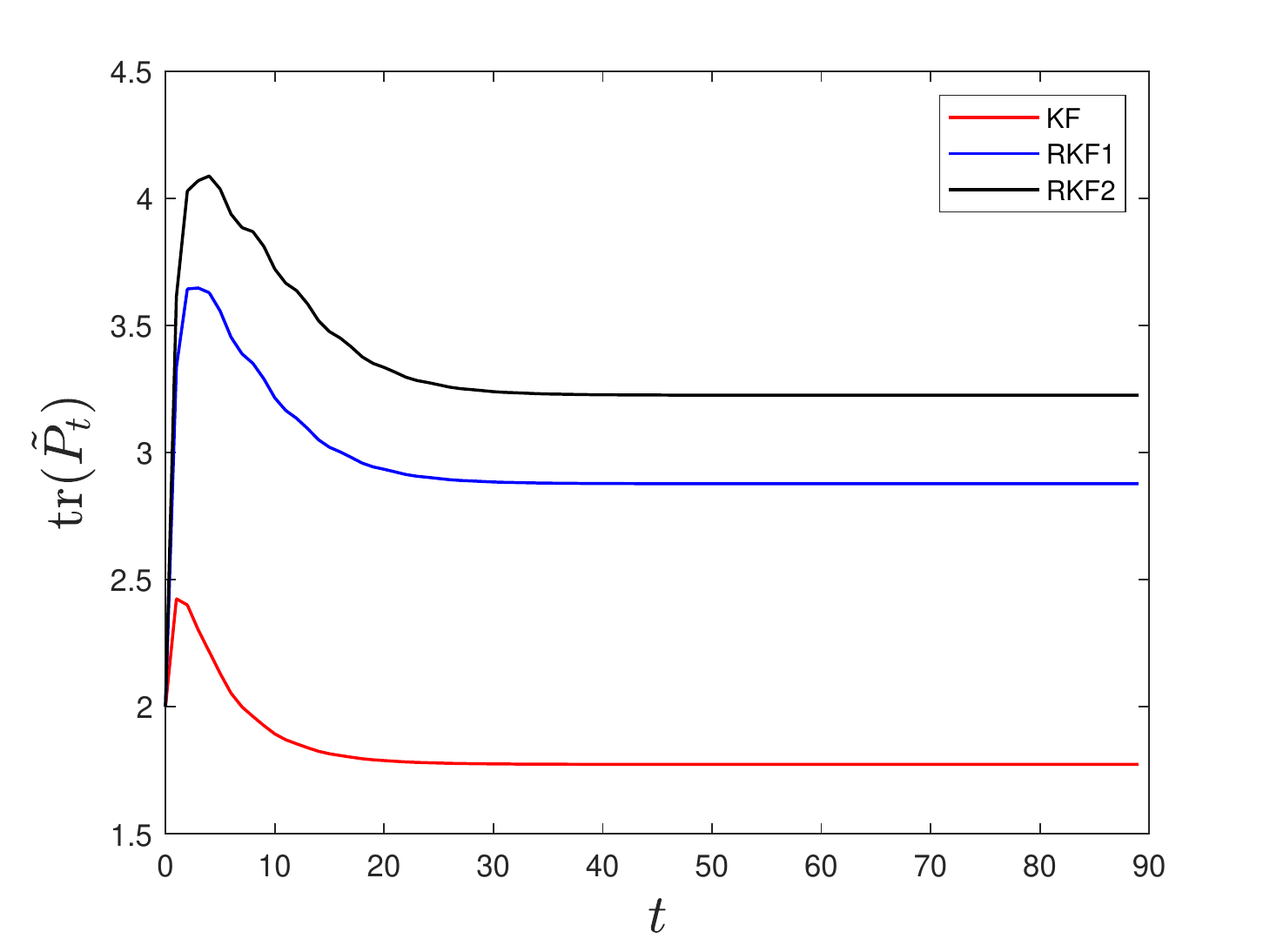}
\caption{Trace of $\tilde P_t$ for KF, RKF1 with $c=5 \cdot 10^{-2}$, and RKF2 with $c=8 \cdot 10^{-2}$. Recall that $\tilde P_t=P_t$ for KF.} \label{fig1}
\end{figure} \begin{figure}[h]
\centering
\includegraphics[width=0.5\textwidth]{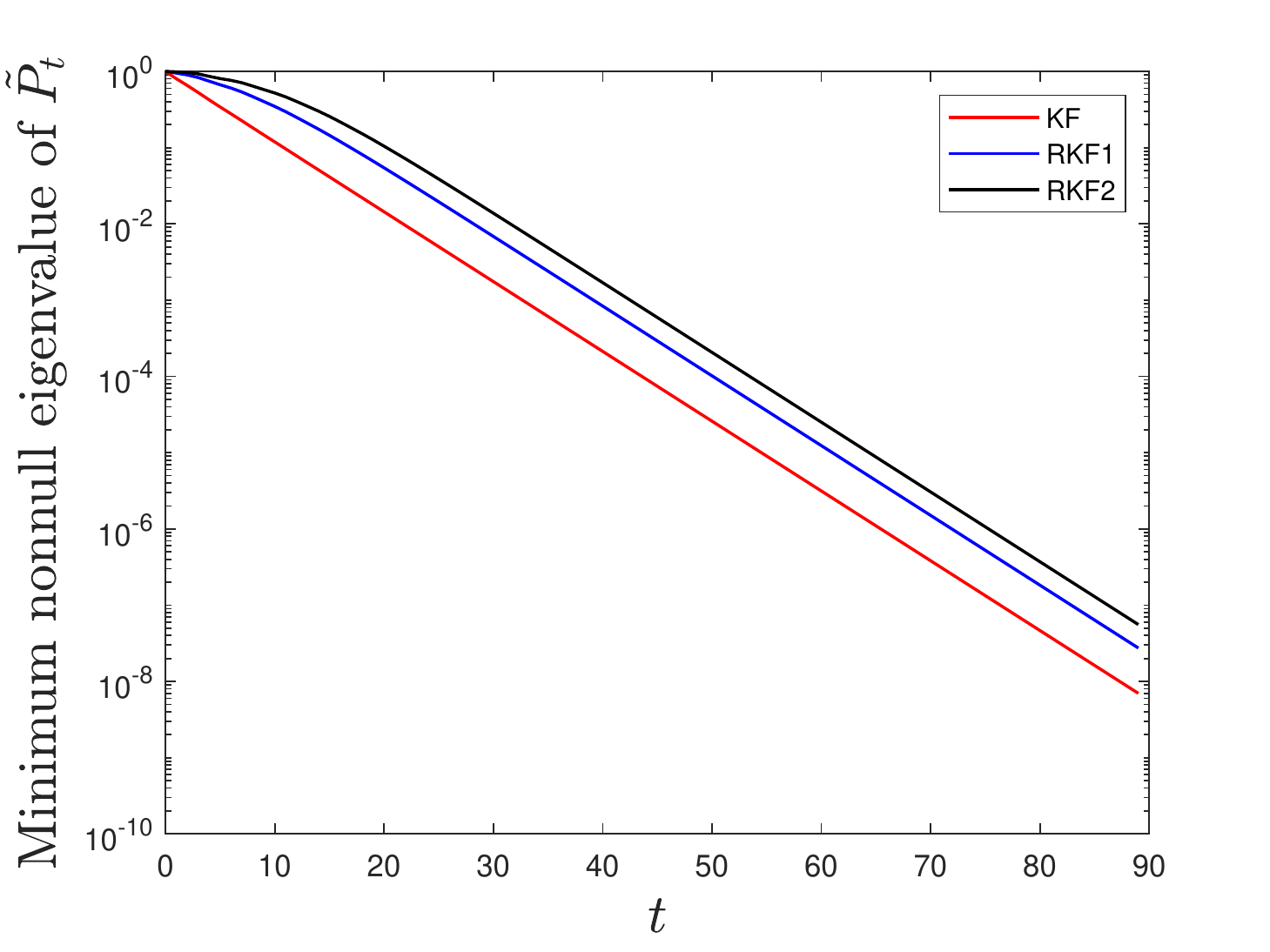}
\caption{Minimum nonnull eigenvalue of $\tilde P_t$ as a function of $t$ for KF, RKF1 with $c=5 \cdot10^{-2}$, and RKF2 with $c=8 \cdot 10^{-2}$.} \label{fig2}
\end{figure}
\begin{figure}[h]
\centering
\includegraphics[width=0.5\textwidth]{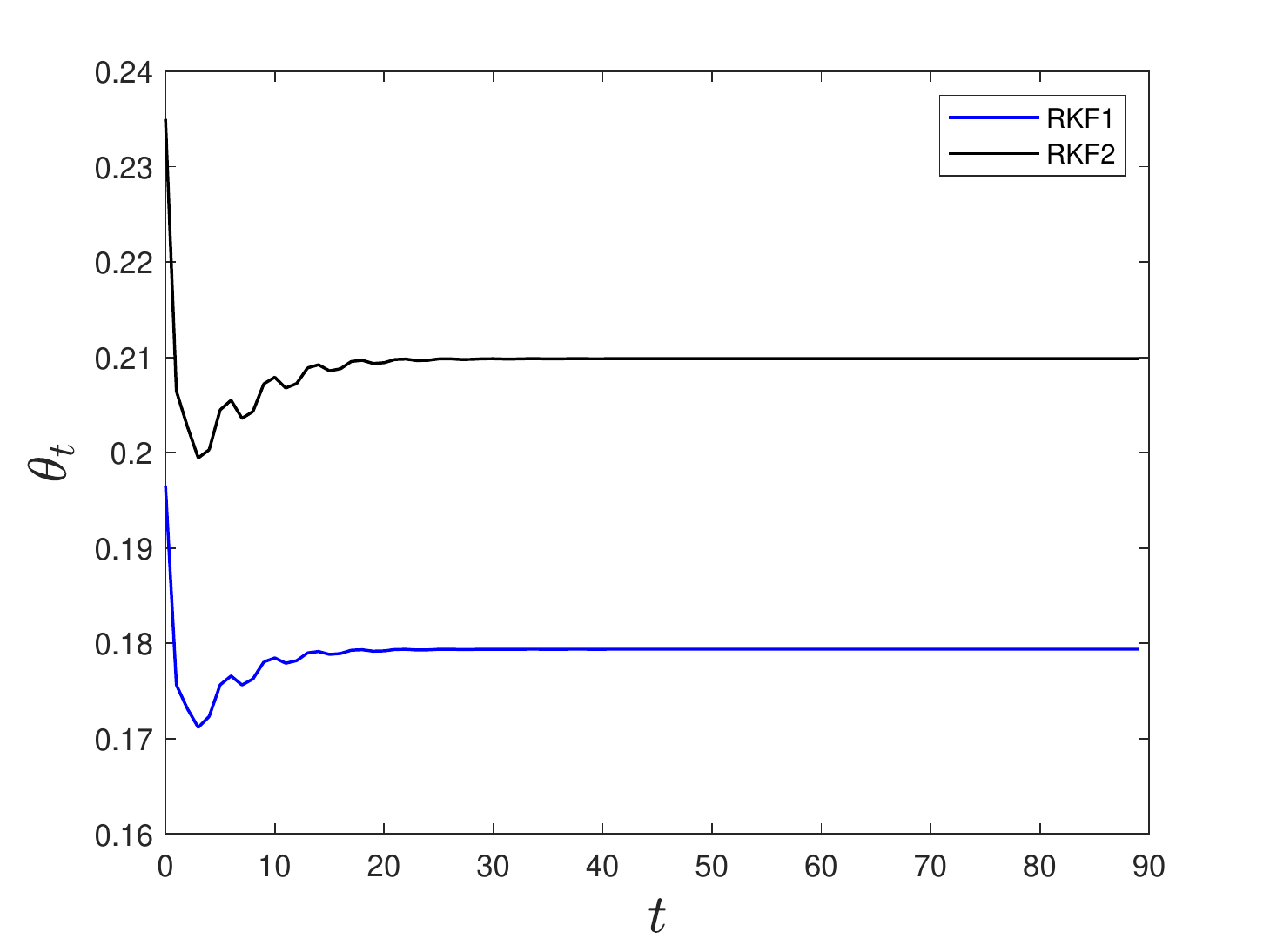}
\caption{Risk sensitivity parameter $\theta_t$ as a function of $t$ for RKF1 with $c=5 \cdot 10^{-2}$, and RKF2 with $c=8 \cdot 10^{-2}$.} \label{fig6}
\end{figure}  Finally,  Figure \ref{fig6} shows the risk sensitivity parameter $\theta_t$ as a function of $t$. As expected, it converges to a constant value and the larger $c$ is, the larger $\theta_t$ is.

Next, we evaluate the performance of RKF2 and KF assuming that the actual model corresponds to (\ref{m_nomi}). Let $e_{t+1}:=x_{t+1}-\hat x_{t+1}$ be the prediction error using a predictor of the form
\begin{align} \label{pred_arb}\hat x_{t+1}=A\hat x_t+G^\prime_t(y_t-C\hat x_t ). \end{align} Then, it is not difficult to see that $e_{t+1}=(A-G_t^\prime C)e_t+(B-G_t^\prime D)v_t$. Thus, $e_t$ is a process with zero mean and with covariance matrix $V_t:=\mathbb E[e_te_t^T]$. The latter is given by solving the following Lyapunov equation
$$ V_{t+1}=(A-G^\prime_t C) V_t(A-G^\prime_t C)^T+(B-G_t^\prime D) (B-G^\prime_t D)^T.$$
\begin{figure}
\centering
\includegraphics[width=0.5\textwidth]{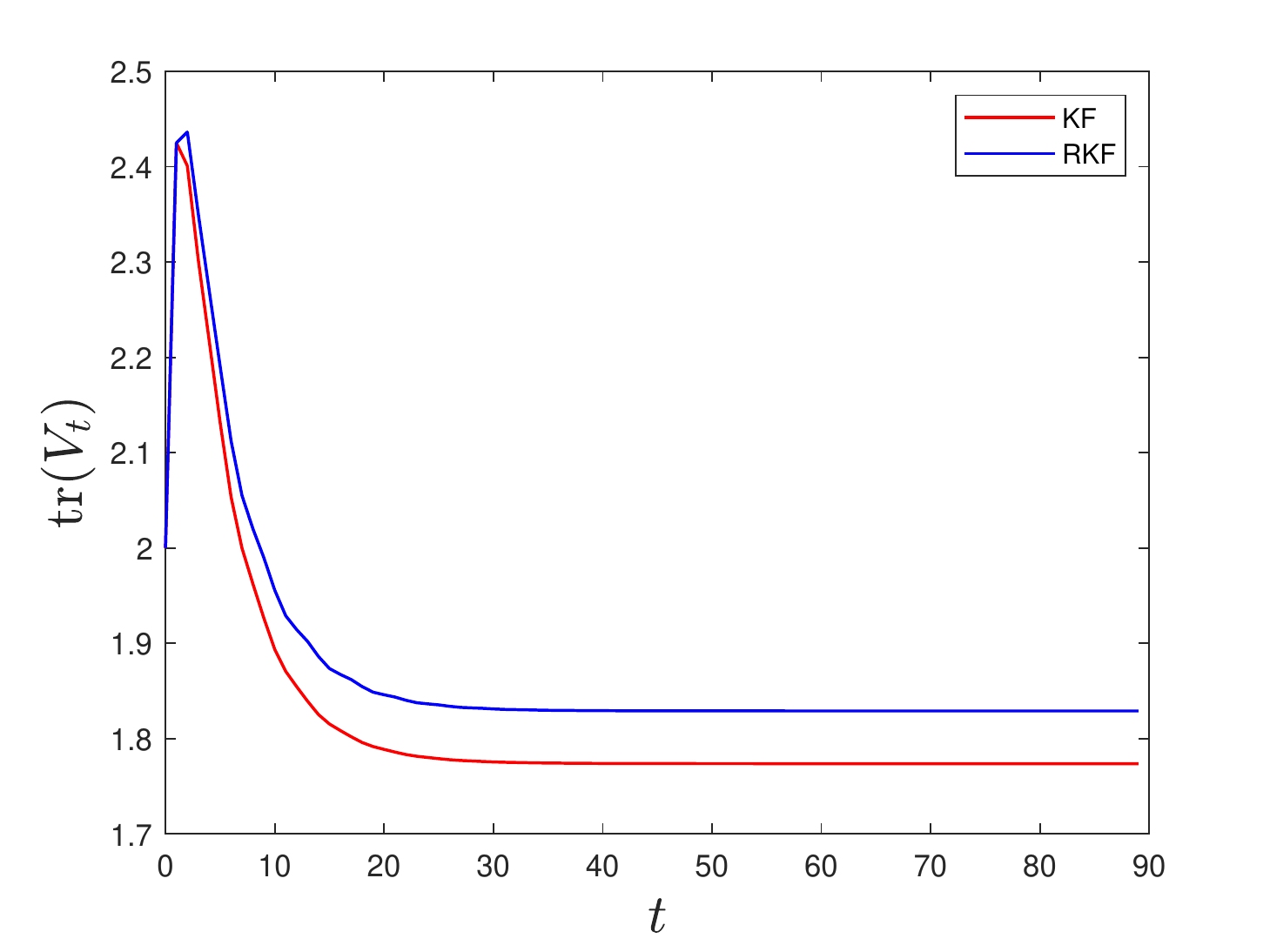}\caption{Scalar variance of the prediction error for KF and RKF2 with $c=8 \cdot 10^{-2}$ under the assumption that the actual model coincides with the nominal model.} \label{fig3a}
\end{figure}Figure \ref{fig3a} shows $\tr(V_t)$ for RKF2 and KF. As expected, KF performs better than RKF2. Indeed, the former has been designed to be optimal for (\ref{m_nomi}).

Finally, we compare the performance of RKF2 and KF  using the actual model
\begin{equation}\begin{aligned}\label{m_lf}
\xi_{t+1} &=\tilde{A} \xi_{t}+\tilde{B} \varepsilon_{t} \\
y_{t} &=\tilde{C} \xi_{t}+\tilde{D} \varepsilon_{t}
\end{aligned}\end{equation}
where
$$\begin{aligned}
&\tilde{A}
=\left[\begin{array}{cccccc}
\gr 1.1 & \gr  0.1 & \gr  0.1 & \bl 0.1033 & \bl 0.0663 & \bl 0.0295 \\
\gr  0 & \gr 0.6364 & \gr -0.6364 &  0 &  0 & 0 \\
0\gr  & \gr 0.6364 & \gr  0.6364 &  0 &  0 & 0\\
0 & 0 & 0 &  0.4365 & 0.2131 & 0.1503\\
0 & 0 & 0 &  0 & 0.6364 & -0.6364\\
0 & 0 & 0 &  0 & 0.6364 & 0.6364\\
\end{array}\right]\\
&\tilde{B}=\left[\begin{array}{cc}
\rd 1.182  & \gr0 \\
\gr 0  &\gr 0\\
\gr 0 & \gr0\\
1.4188 &  -0.92\\
0&0\\
0&0
\end{array}\right]\\
&\tilde{C}=\left[\begin{array}{llllll}
\gr 1 & \gr 0 & \gr 0 & \bl -0.0868 & \bl -0.0557 & \bl -0.0247
\end{array}\right]\\
&\tilde{D}=\left[\begin{array}{llll}
\rd -0.2821  & \rd 1.0956
\end{array}\right],\\
\end{aligned}$$ $\xi_t=\left[\begin{array}{ll} x_t^T & \eta_t^T \end{array} \right]^T$, and $\varepsilon_t$ is normalized white Gaussian noise. Note that $x_t$ is the actual state, {while $\eta_t$ is a perturbation process. Model (\ref{m_lf}) is a perturbed version of (\ref{m_nomi}). Indeed, the green parts in $\tilde A$, $\tilde B$, $\tilde C$, $\tilde D$ correspond to the non-perturbed parts of matrices $A$, $B$, $C$ and $D$; the red parts in $\tilde A$, $\tilde B$, $\tilde C$, $\tilde D$ correspond to the perturbed parts of matrices $A$, $B$, $C$ and $D$; the blue parts are the terms coupling $x_t$ and $y_t$ with the perturbation process $\eta_t$. Consider a predictor of the form (\ref{pred_arb}). Let $e_{t}=x_t-\hat x_t$ denote the prediction error of such a predictor under the actual model in (\ref{m_lf}). Since the submatrices in $\tilde A$ and $\tilde C$ corresponding to $A$ and $C$ are not perturbed, then it is not difficult to see that
$$ \tilde \xi_{t+1}=(\tilde A- \tilde G_t^\prime \tilde C )\tilde \xi_t+(\tilde B-\tilde G_t^\prime \tilde D)\varepsilon_t$$
where $\xi_t:=[\, e_t^T\; \eta_t^T\,]^T$, $\tilde G_t^\prime=[\, (G_t^\prime)^T \; 0\,]^T$. Thus, $e_t$ is a process with zero mean and with covariance matrix $V_t:=\mathbb E[e_te_t^T]$. The latter is given by solving the following Lyapunov equation
$$ \Pi_{t+1}=(\tilde A-\tilde  G^\prime_t \tilde C) \Pi_t(\tilde A-\tilde G^\prime_t \tilde C)^T+(\tilde B-\tilde G_t^\prime\tilde  D ) (\tilde B-\tilde G^\prime_t \tilde D)^T$$
and $$ \Pi_t= \left[\begin{array}{cc} V_t & \star  \\ \star & \star  \end{array}\right].$$
\begin{figure}
\centering
 \includegraphics[width=0.5\textwidth]{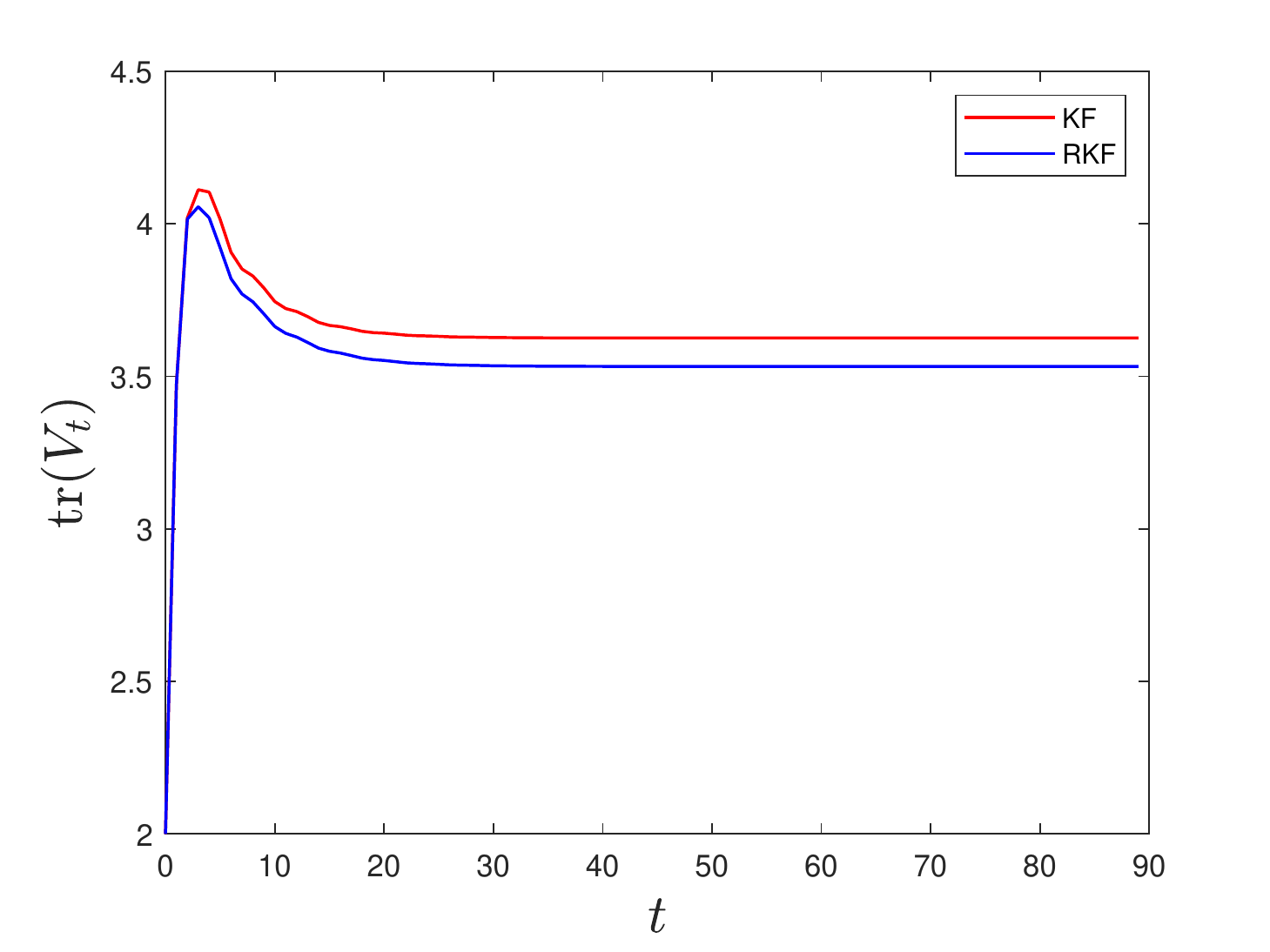}
\caption{Scalar variance of the prediction error  for KF and RKF2 with $c=8 \cdot 10^{-2}$ under the assumption that the actual model is different from the nominal one.} \label{fig3}
\end{figure}Figure \ref{fig3} shows $\tr(V_t)$ for RKF2 and KF. In this case the former performs better than the latter.}

\section{Conclusions}\label{sec_5}
{ In this paper, we consider a robust static estimation problem in the case that the nominal density is Gaussian and possibly degenerate. We apply such a result to design a robust Kalman filter which can be used also in the case that the reachability assumption does not hold. Of course, there are many aspects that have not taken into account yet. More precisely, in our research agenda now there are the following questions that we are addressing:
 \begin{itemize}
 \item The least favorable model has been assumed to be Gaussian. However, in view of the results in \cite{robustleastsquaresestimation}, we believe our conclusions also holds in the case that the ambiguity set contains non-Gaussian probability densities.
 \item The minimax problem also provides the least favorable density $\tilde f^0_t$. Under the assumption that $f_t$ and $\tilde f_t$ are non-degenerate, Levy \& Nikoukhah showed that it is possible to characterize the corresponding dynamic model in a finite simulation horizon, \cite{ROBUST_STATE_SPACE_LEVY_NIKOUKHAH_2013}. Our conjecture is that it is possible to adapt these arguments to our case.
 \item The simulation results show that the robust Kalman filter seems to converge in the case of constant parameters. Drawing inspiration from the non-degenerate case in \cite{ZORZI_CONTRACTION_CDC,CONVTAU}, we believe it is possible to prove that the risk sensitive-like Riccati iteration in Algorithm \ref{code:recentEnd} (steps \ref{RIcc_step} and \ref{RIcc_step2}) converges provided that the tolerance is sufficiently small and the system is stabilizable and detectable.
 \end{itemize}}



\end{document}